\documentclass{amsart}
\usepackage{graphicx}
\usepackage{amsmath}
\usepackage{amssymb}
\usepackage{amsthm}
\usepackage{wrapfig}

\newtheorem{thm}{Theorem}[section]
\newtheorem{cor}[thm]{Corollary}
\newtheorem{conj}[thm]{Conjecture}
\newtheorem{prob}[thm]{Problem}
\newtheorem{lem}[thm]{Lemma}

\theoremstyle{definition}
\newtheorem{defn}[thm]{Definition}
\theoremstyle{remark}

\numberwithin{equation}{section}

\begin{document}

\title[Density of Iterated Line Intersections]{A Result About the Density of Iterated \\ Line Intersections in the Plane}%
\author{Christopher J. Hillar}%
\address{Department of Mathematics, University of California, Berkeley, CA 94720.}
\email{chillar@math.tamu.edu}

\author{Darren L. Rhea}%
\address{Department of Mathematics, University of California, Berkeley, CA 94720.}
\email{drhea@math.berkeley.edu}

\thanks{The work of the first author is supported under a National
Science Foundation Graduate Research Fellowship.}


\begin{abstract}
Let $S$ be a finite set of points in the plane and let
$\mathcal{T}(S)$ be the set of intersection points between pairs
of lines passing through any two points in $S$.  We characterize
all configurations of points $S$ such that iteration of the above
operation produces a dense set.  We also discuss partial results
on the characterization of those finite point-sets with rational
coordinates that generate all of $\mathbb Q^2$ through iteration
of $\mathcal{T}$.
\end{abstract}
 \maketitle
\section{Introduction}

Let $S$ be a set of points in the plane and let $L = \{L_{i}\}_{i
\in I}$ be the set of lines between pairs of points in $S$.
Consider the following operation on $S$:

\begin{equation}
\mathcal{T}(S) = \bigcup\limits_{i \ne j} {L_i  \cap L_j } \subseteq
\mathbb R^2.
\end{equation}
In other words, $\mathcal{T}(S)$ is the set of intersection points
between pairs of distinct lines in $L$.  If $S$ consists of $n$
collinear points (or no points at all), then the union above is
empty; so to keep the notation consistent, we set $\mathcal{T}(S)
= \emptyset$ for these cases.

\begin{figure}[!htbp]
\begin{center}
\includegraphics[scale=0.45]{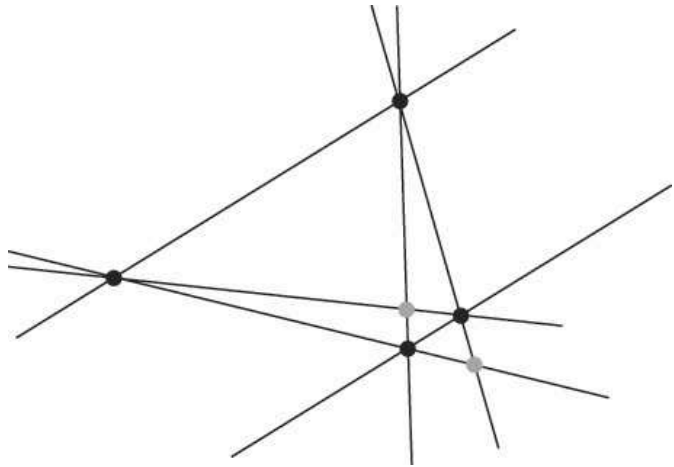}
\end{center}
\caption{$\mathcal{T}(S)$ for a set of points $S$ that form a trapezoid.}
\label{fig.trapezoid}
\end{figure}

As a simple example of the operation $\mathcal{T}$, let $S$ 
consist of four black points that are the vertices of a trapezoid 
as in  Figure \ref{fig.trapezoid}.  Then, $\mathcal{T}$$(S)$ 
consists of the original four points along with two additional 
ones shown in gray.  It should be clear that for a set of points not
all collinear, we have $S \subseteq \mathcal{T}(S)$. Moreover,
$\mathcal{T}(S)$ is finite for finite sets $S$.  We are
interested here in the iterations, $\mathcal{T}^i(S)$, and
specifically, the limiting behavior of such operations on
arbitrary finite sets $S$. The study of such phenomenon naturally
leads to the notion of the order of a set $S$, which we define
below.  As a matter of convention, we set $\mathcal{T}^0(S)$ = S. 

\begin{defn}
Let $S$ be a set of points in $\mathbb R^2$.  The \emph{order} of
$S$ is the smallest positive integer $n$ such that
$\mathcal{T}^n(S) = \mathcal{T}^{n-1}(S)$.  If there is no such $n$, then the order
of $S$ is defined to be $\infty$.
\end{defn}

For example, the order for a set of points forming the vertices of a square is 2.  If
the order of a set $S$ is 1, then we call $S$ \emph{fixed} under
$\mathcal{T}$.  A set $S$, therefore, has finite order if and only if
$\mathcal{T}^n(S)$ is fixed for some nonnegative integer $n$.

\begin{prob}\label{fixedprob}
Describe the finite point-sets that have finite order.
\end{prob}

Before discussing the answer to this problem
(in Section \ref{fixedsetssec}), we describe a nontrivial infinite
point-set that has finite order.  Let $S$ be the set of rational
points on the unit circle, $x^2 + y^2 = 1$.  For a given $P \in
\mathbb Q^2$, choose two points $A$ and $B$ in $S$ such that $PA$
and $PB$ are not tangent to the unit circle.  Then, if $C$ and $D$
are the points of intersection of $PA$ and $PB$ (respectively)
with the circle, it turns out \cite[p. 249]{numtheory} that $C$
and $D$ are both rational. It follows that $P \in \mathcal{T}(S)$
for every $P \in \mathbb Q^2$, and thus \[\mathcal{T}^2(S) =
\mathcal{T}(\mathbb Q^2) = \mathbb Q^2 = \mathcal{T}(S).\]

Excluding the sets of finite order, it follows that iteration of
$\mathcal{T}$ produces a strictly increasing chain of sets of
points in the plane.  In light of this observation, a natural
question is whether we arrive at a dense set of points by such a
procedure. In other words, is $\bigcup\nolimits_{i \geq 0}
{\mathcal{T}^{i} \left( S \right)}$ dense in $\mathbb R^2$? A more
difficult but related question is whether we get all of $\mathbb
Q^2$ when $S$ consists of only rational points.  We address both
of these questions with a complete answer to the first in Section
\ref{denseproof} and some partial results for the second in Section 
\ref{rationalsection}.

\begin{thm}\label{densethm}
Let $S$ be a finite set of points in the plane.  Then, $S$ has
infinite order if and only if $\bigcup\nolimits_{i \geq 0}
{\mathcal{T}^{i} \left( S \right)}$ is dense in $\mathbb R^2$.
\end{thm}

The answer to Problem \ref{fixedprob} found in Corollary
\ref{fixedthmcor} below, therefore, gives a complete
characterization of when iterated line intersections are dense.

\begin{cor}\label{fulldensechar}
Let $S$ be a finite set of points in the plane.  Then, 
$\bigcup\nolimits_{i \geq 0} {\mathcal{T}^{i} \left( S \right)}$ 
is dense in $\mathbb R^2$ if and only if $S$ is not one of the 
following sets:
\begin{enumerate}
    \item The empty set.
    \item A set of collinear points.
    \item A set of collinear points with one additional noncollinear point.
    \item The vertices of a parallelogram.
    \item The vertices of a parallelogram and the intersection of its two
    diagonals.
\end{enumerate}
\end{cor}

In the rational case, we conjecture a more exact result.

\begin{conj}\label{rationalthm}
Let $S$ be a finite set of points in the plane with rational
coordinates. Then, $S$ has infinite order if and only if
$\bigcup\nolimits_{i \geq 0} {\mathcal{T}^{i} \left( S \right)} =
\mathbb Q^2.$
\end{conj}

As a step in the direction of this conjecture,
we offer the following; its proof can be found in Section \ref{rationalsection}.

\begin{thm}\label{rationalparallelthm}
Let $R, P, Q, T \in S$ be rational points in the plane with $RQ$
and $PT$ parallel and suppose that $RP$ is not parallel to $QT$.
Then, $\bigcup\nolimits_{i \geq 0} {\mathcal{T}^{i} \left( S
\right)} = \mathbb Q^2$.
\end{thm}

Though we were not motivated by any other particular work, we should remark that a 
similar question posed by Fejes-Toth (with circles replacing lines) was addressed by 
Bezdek and Pach in \cite{pach}, and related results can 
also be found in the papers \cite{barany, king}.  Additionally,
Theorem \ref{densethm} has also been 
discovered recently (independently) by Ismailescu and Radoicic \cite{otherproof}.

\section{Finite Fixed Sets}\label{fixedsetssec}

We begin by characterizing sets of finite order.  Although one may deduce 
the main result of this section from Lemmas \ref{trianglelemma} and 
\ref{trianglelemma2} in Section \ref{denseproof}, the methods employed here are less 
cumbersome and might be of independent interest.  We will need 
the following result from elementary geometry.

\begin{thm}[The Sylvester-Gallai Theorem]\label{sylthm}
For every set of $n$ noncollinear points in the plane, there
exists a line that contains exactly two of the points.
\end{thm}

Although this fact seems intuitively obvious, its proof eluded 
even Sylvester, and it was only solved (in published form) some $50$ 
years after being posed by him \cite{solvedsylvester}.  We refer the reader to \cite{bookproof} for
more details.  We are ready to approach Problem \ref{fixedprob}.

\begin{thm}\label{fixedthm}
A finite set $S$ fixed under $\mathcal{T}$ must be one of the following configurations:
\begin{enumerate}
    \item The empty set.
    \item A set of collinear points with one additional noncollinear point.
    \item The vertices of a parallelogram and the intersection of its two
    diagonals.
\end{enumerate}
\end{thm}
\begin{proof}
Let $S$ be a set of $n$ noncollinear points in the plane that is
fixed by $\mathcal{T}$.  Using Theorem \ref{sylthm}, there exists a
line intersecting $S$ in exactly two points $P$ and $Q$.  By
assumption, there is some other point $X$ not on this line, and
we can choose $X$ so that its altitude from $PQ$ is largest.  If
all other points lie on the line $XP$ or if all of them lie on
$XQ$, then we are in configuration (2) above.  The remaining
possibilities break up into two cases.

\begin{figure}[!htbp]
\begin{center}
\includegraphics[scale=0.45]{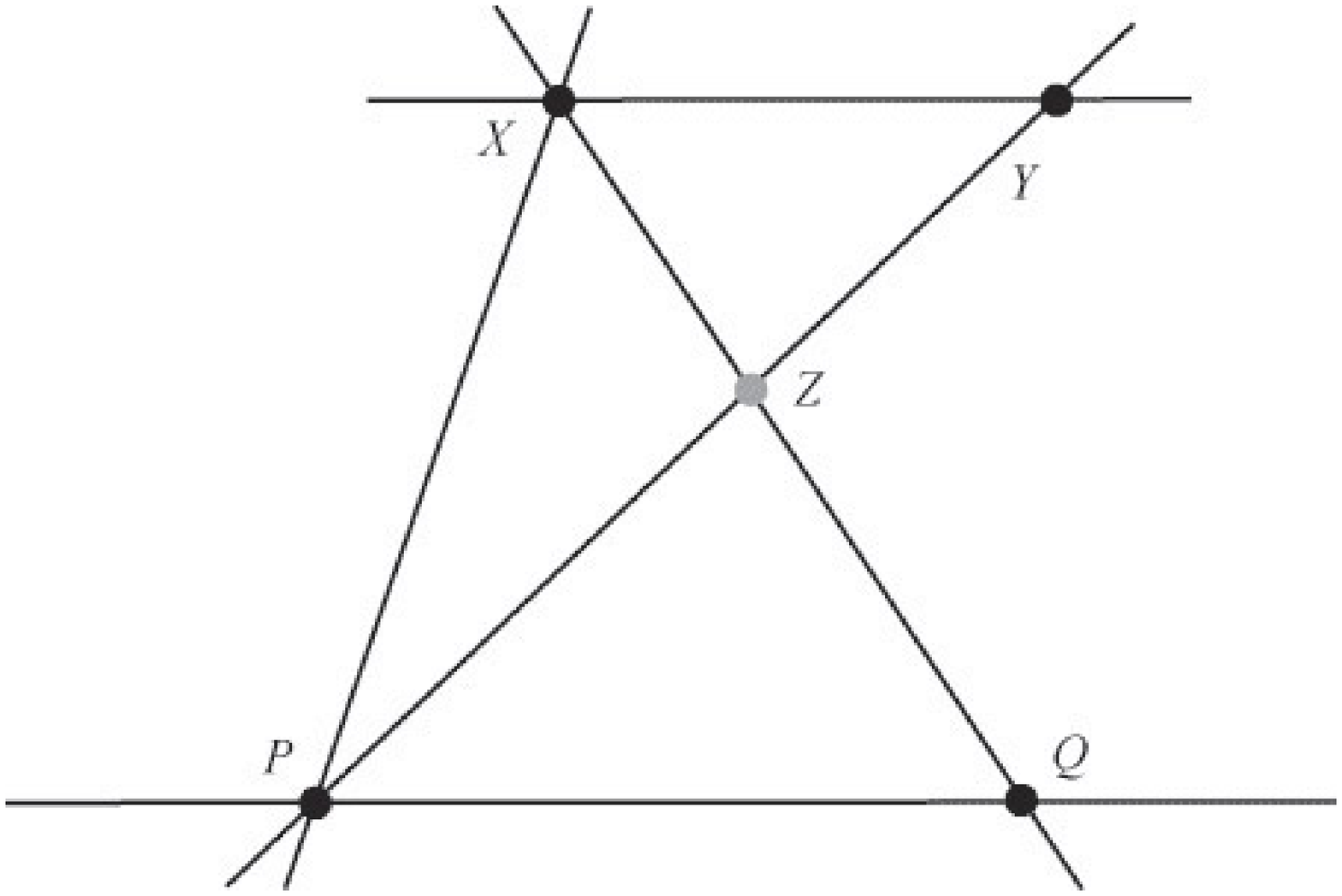}
\end{center}
\caption{Case 1 in the proof of Theorem \ref{fixedthm}} \label{fig.fixedproof1}
\end{figure}

\underline{Case 1}:  There is a point $Y \in S$ not on $XP$ and
not on $XQ$.


We first claim that $Y$ must lie on the line through $X$ that is parallel 
to $PQ$. Indeed, any other
position for $Y$ would give rise to an intersection between $XY$
and $PQ$ that is not $P$ or $Q$, contrary to our use of Theorem
\ref{sylthm} and our assumption that $\mathcal{T}(S) = S$.
Relabeling if necessary, Figure \ref{fig.fixedproof1} depicts the situation.
Since $S$ is fixed, the intersection point, $Z$, of $XQ$ and $PY$ is in
$S$. It follows that $XP$ and $YQ$ 
must be parallel (otherwise, if
$W$ is the intersection point of $XP$ and $YQ$, then $ZW$ would
intersect $PQ$). Finally, it is easy to see that there can be no other points
in $S$ by our choice of $P$ and $Q$.

\begin{figure}[!htbp]
\begin{center}
\includegraphics[scale=0.45]{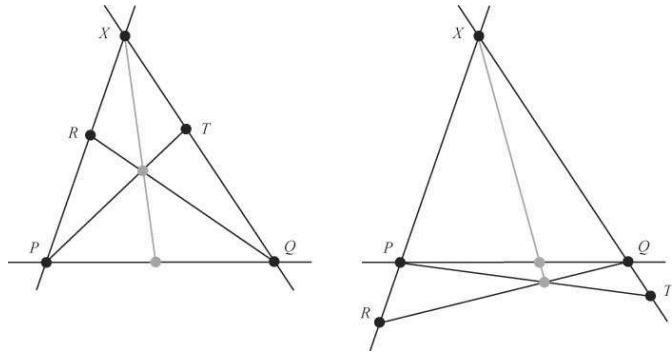}
\end{center}
\caption{Case 2 in the proof of Theorem \ref{fixedthm}} \label{fig.fixedproof2}
\end{figure}

\underline{Case 2}: Every point in $S$ lies on one of the lines
$XP$ or $XQ$.

If $S$ is not a configuration of type (2), then there are points $R,T \in S$ such
that $R$ is on the line $XP$, $T$ is on the line $XQ$, and $R,T$
are not $X,P$, or $Q$.  By the assumption on $X$ and the line $PQ$, 
only two configurations for $R$ and $T$ are possible; these are 
depicted in Figure \ref{fig.fixedproof2}.  In both cases, two iterations
of $\mathcal{T}$ give rise to a point in $S$ on the line $PQ$, a contradiction.
Therefore, no fixed point-sets other than those of configuration (2) may 
take this form.  This completes the proof.
\end{proof}

\begin{cor}\label{fixedthmcor}
The finite point-sets with finite order are
\begin{enumerate}
    \item The empty set.
    \item A set of collinear points.
    \item A set of collinear points with one additional noncollinear point.
    \item The vertices of a parallelogram.
    \item The vertices of a parallelogram and the intersection of its two
    diagonals.
\end{enumerate}
\end{cor}

\begin{proof}
Let $S$ be a finite set in $\mathbb R^2$ with order $n$.
Applying Theorem \ref{fixedthm}, it
follows that $R = \mathcal{T}^{n-1}(S)$ must be one of three types.  
When $R$ is empty, then $S$ is either itself empty or a set of 
collinear points.  Similarly, a set $R$ of collinear points with one
additional point can only be obtained from a set $S$ that is the same as $R$. 
Finally, when $R$ forms a parallelogram with the intersection of its diagonals, 
the set $S$ must either be $R$ or $R$ without its diagonal intersection.
\end{proof}

\section{The Density Theorem}\label{denseproof}

Before proving Theorem \ref{densethm}, we record the following
technical lemmas, the first of which provides a useful
characterization of sets of infinite order.   For ease of 
presentation, we say that a point is \emph{strictly contained} in a 
set $K$ if it is located in its interior.

\begin{lem}\label{trianglelemma}
Let $S$ be a finite set of infinite order.  Then, there exists $n
\in \mathbb N$ such that $\mathcal{T}^n(S)$ contains a subset of 4
points in which 3 of the points are noncollinear and
the fourth point is strictly contained in the triangle determined
by these 3 points.
\end{lem}

\begin{proof}
We consider the number of vertices $v$ on the convex hull $H$
of $S$.  When $v = 2$, the set $S$ cannot have 
infinite order.  So suppose that $v = 3$.  If there is a point of $S$ 
strictly contained inside $H$, then we are done.
Otherwise, since $S$ has infinite order, there must be two points of $S$
on different edges of $H$.  An iteration of $\mathcal{T}$ then produces our
desired point.

\begin{figure}[!htbp]
\begin{center}
\includegraphics[scale=0.45]{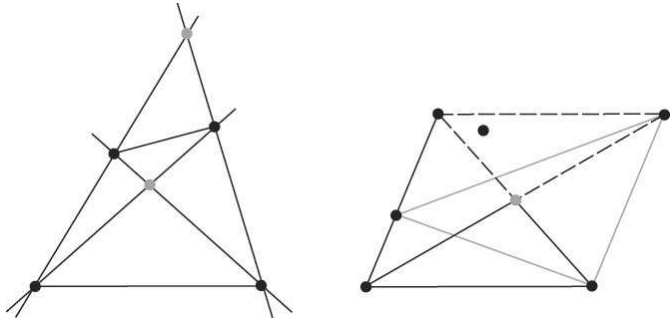}
\end{center}
\caption{Four vertices on the convex hull of $S$} \label{fig.trianglempic}
\end{figure}

Assume now that $H$ has exactly four vertices.  If these vertices do
not form a parallelogram, then one iteration of $\mathcal{T}$ gives us
what we want (see Figure \ref{fig.trianglempic}).  Otherwise, there is
a point in $S$ which is not a vertex of $H$ and not the intersection
of the diagonals of the quadrilateral determined by $H$.  Again in this 
case, one iteration of $\mathcal{T}$ (giving us the intersection
of the two diagonals of $H$) produces the desired result.

Finally, if $v > 4$, then we proceed as follows.
Pick two adjacent vertices $A$ and $B$.  There must be two other 
vertices $C$ and $D$ such that the edges $AB$ and $CD$ are 
not parallel ($H$ has at least $5$ vertices and is convex).  This
reduces the problem to the case of $4$ vertices not forming 
a parallelogram (encountered above) and completes the proof of the lemma.
\end{proof}

Our next result allows one to produce a convergent, nested sequence of triangles.

\begin{lem}\label{trianglelemma2}
Let $A$, $B$, and $C$ be noncollinear points, and let $P$
be a point strictly inside $\triangle ABC$.  Then, there
exist triangles $\triangle A_nB_nC_n$ ($n =
1,2,\ldots$) strictly containing $P$ such that $\mathop {\lim
}\limits_{n \to \infty } A_n  = \mathop {\lim }\limits_{n \to
\infty } B_n = \mathop {\lim }\limits_{n \to \infty } C_n = P$, and for each $n$,
\[A_n,B_n,C_n \in \bigcup\limits_{j = 0}^\infty {
\mathcal{T}^{(j)} \left(\{A,B,C,P\}  \right )}.\]
\end{lem}

\begin{figure}[!htbp]
\begin{center}
\includegraphics[scale=0.45]{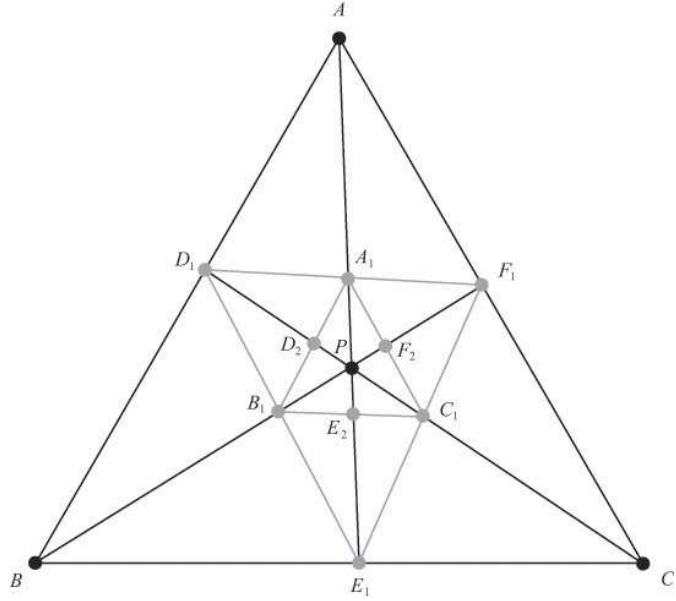}
\end{center}
\caption{Nested triangle iteration} \label{fig. nestedtriangles}
\end{figure}

\begin{proof}
Given a triangle $\triangle ABC$ and a point $P$ strictly contained in it, 
we may construct the vertices of another triangle containing this 
point by intersecting the lines $AP$, $BP$, and $CP$ 
with the edges of $\triangle ABC$.
Iterating this procedure produces a nested sequence of
triangles strictly containing $P$ with vertices in
$\bigcup\nolimits_{j \geq 0} {\mathcal{T}^{(j)} \left( \{A,B,C,P
\}\right)}$ (see Figure \ref{fig. nestedtriangles}).  This sequence 
contains two types of triangles; we label the odd iterates $\triangle D_nE_nF_n$, 
while even iterates are denoted by $\triangle A_nB_nC_n$.
Here, the $A_n$ (resp. $B_n$, $C_n$) are labeled so that they 
are the ones on the line $AP$ (resp. $BP$, $CP$).
We claim the vertices of the triangles $\triangle A_nB_nC_n$ 
all converge to $P$.  

\begin{figure}[!htbp]
\begin{center}
\includegraphics[scale=0.5]{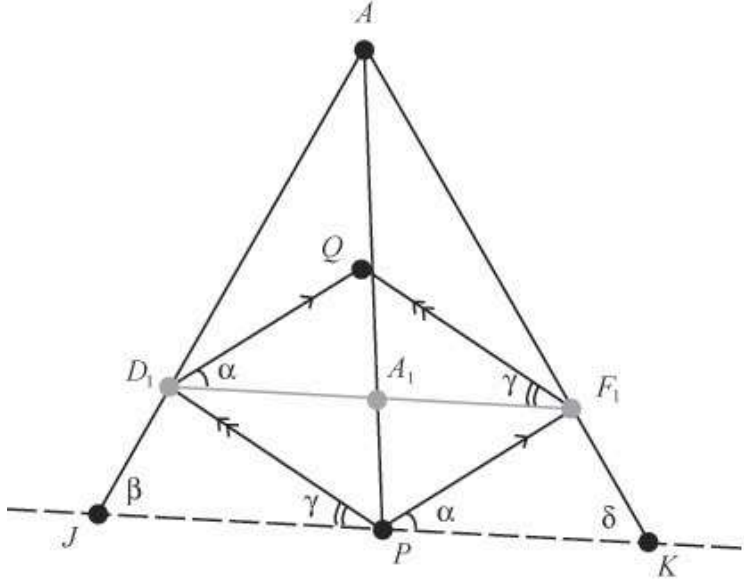}
\end{center}
\caption{Iterations decrease triangle areas} \label{fig. nestedtriangles2}
\end{figure}

To verify this assertion, it suffices to show that $|A_1P| < |AA_1|$, $|B_1P| < |BB_1|$, 
and $|C_1P| < |CC_1|$.  Without loss of generality, we prove that $|A_1P| < |AA_1|$.
Reducing further, we observe that it is enough to show that the area of
$\triangle PD_1F_1$ is less than the area of $\triangle AD_1F_1$ 
(drop altitudes to $D_1F_1$ from $A$, $P$ and compare similar triangles).
Next, draw the line $JK$ that is parallel to $D_1F_1$ and passes through $P$, 
and label the angles formed as in Figure \ref{fig. nestedtriangles2}.  
Since $F_1P$ and 
$AJ$ (resp. $D_1P$ and $AK$) intersect at $B$ (resp. $C$), it follows that 
$\alpha < \beta$ and $\gamma < \delta$.  Therefore, when we form the
triangle $\triangle QD_1F_1$ that is congruent to $\triangle PD_1F_1$, 
it must lie entirely inside $\triangle AD_1F_1$.  This finishes the proof.
\end{proof}

\begin{lem}\label{densetrianglelem}
Let $A$, $B$, $C$ be noncollinear points in the plane.  If $K$ is a dense set
of points in $\triangle ABC$, then $\mathcal{T}(K)$ is
a dense set of points in the entire plane.
\end{lem}
\begin{proof}
Let $P$ be a point in the plane, and let $Q_1, Q_2$ and $R_1, R_2$ be points 
strictly inside $\triangle ABC$ such that $Q_1Q_2$ and $R_1R_2$
intersect at $P$.  Since $K$ is dense in $\triangle ABC$,
there are a sequence of points $Q_{1n},Q_{2n} \in K$ and $R_{1n}, R_{2n} \in K$ 
that converge to $Q_1, Q_2$ and $R_1, R_2$, respectively.  Since the intersection
of two lines formed by four points is continuous in the four points (the intersection is
a rational function in the coordinates of the four points), it follows that the intersections
of $Q_{1n}Q_{2n}$ and $R_{1n}R_{2n}$ (which are in $\mathcal{T}(K)$) converge to
$P$.  This completes the proof.
\end{proof}

We are ready to prove Theorem \ref{densethm}.

\begin{proof}[Proof of Theorem \ref{densethm}]
The if-direction ($\Leftarrow$) in the theorem statement is immediate.  Therefore, let $S$ be
a finite set of infinite order.  Using Lemma \ref{trianglelemma},
there exists $n \in \mathbb N$ such that $\mathcal{T}^n(S)$
contains a triangle of vertices and a fourth point strictly
contained in the triangle determined by these 3 vertices.  We
claim that iteration of $\mathcal{T}$ on these 4 points produces a
dense set of points in the triangle.  The theorem then follows from Lemma \ref{densetrianglelem}.

Let $A$, $B$, and $C$ be the vertices of the triangle strictly
containing $P$.  Suppose that $K = \bigcup\nolimits_{j \geq 0}
{\mathcal{T}^{(j)} \left( A,B,C,P \right)}$ does not contain a
dense set of points in $\triangle ABC$; we will derive a
contradiction.  Using Lemma \ref{trianglelemma2}, we can produce a
sequence of triangles, $\triangle A_iB_iC_i$, with vertices in $K$
such these vertices converge to $P$.  Let $h$ be so
large that the circle centered at $P$ with radius equal to twice
the largest side of $\triangle A_hB_hC_h$ is strictly contained in
$\triangle ABC$.  Since $K$ is not dense in $\triangle ABC$, 
it follows that $K$ cannot be dense in $\triangle A_hB_hC_h$ (again using 
Lemma \ref{densetrianglelem}).

Let $\overline{K}$ be the closure of $K$ and set $W = \overline{K}
\cap \triangle A_hB_hC_h$.  Also, let Int$(\triangle
A_hB_hC_h)$ denote the interior of $\triangle
A_hB_hC_h$.  Since $K$ is not dense in the triangle $\triangle A_hB_hC_h$, 
the (nonempty) open set Int$(\triangle
A_hB_hC_h) \setminus W$ contains an open ball centered at some point $X$
inside $\triangle A_hB_hC_h$.  Consider the set of all closed balls
centered at $X$ that do not intersect $\overline{K}$, and let $r >
0$ denote the supremum over all radii of such balls.  The closed
ball  $\overline{B}(X,r)$ of radius $r$ centered at $X$ must be strictly contained in $\triangle
ABC$ since its interior cannot contain $A_h$, $B_h$, or $C_h$
(they are in $K$) and because of how we chose $h$.

\begin{figure}[!htbp]
\begin{center}
\includegraphics[scale=0.45]{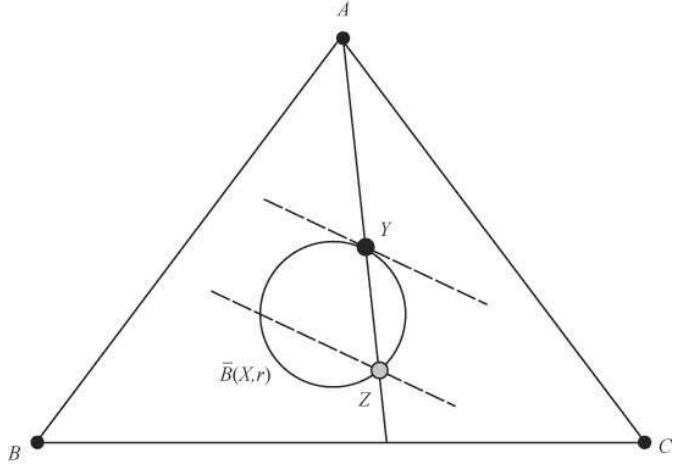}
\end{center}
\caption{Obtaining a contradiction} \label{fig.balltangent}
\end{figure}

By construction of $\overline{B}(X,r)$, there exists a
point $Y \in \overline{K}$ intersecting the boundary of $\overline{B}(X,r)$.
Consider the
lines $AY$, $CY$, and $BY$, and notice that they cannot all be
tangent to the ball $\overline{B}(X,r)$ (there is only one tangent
line through a point on a circle).  Therefore, at least one of
these lines through $Y$, say $AY$, must intersect the interior of
$\overline{B}(X,r)$.  Let $Z$ be the intersection of the line $AY$ 
with the boundary of $\overline{B}(X,r)$ (the point $Z$ need not be in $\overline{K}$). 
The situation is depicted in Figure \ref{fig.balltangent}.  The 
dashed line through $Y$ is the line tangent to the boundary 
of $\overline{B}(X,r)$ at $Y$, while the dashed line through $Z$ is parallel to it.

To continue, we observe the following straightforward fact that was 
discussed in the proof of Lemma \ref{densetrianglelem}: If
$U,V,Q,R \in \overline{K}$ determine two nonparallel lines $UV$
and $QR$, then the intersection point of $UV$ and $QR$ is in
$\overline{K}$.  With this observation in mind, we may use Lemma
\ref{trianglelemma2} to obtain vertices of triangles $\triangle A'_iB'_iC'_i$
in $\overline{K}$ that contain $Y$ and that also converge to $Y$.
None of the vertices $A'_i$, $B'_i$, or $C'_i$ is in the interior of
$\overline{B}(X,r)$ by our choice of $r$. 

Finally, we claim that for large 
enough $n$, the segment $YZ$ must intersect a side of
$\triangle A'_n B'_n C'_n$ in the interior of $\overline{B}(X,r)$,
a contradiction to our assumption on $r$.  To see this, notice that for a large $n$,
at least one of the vertices of $\triangle A'_n B'_n C'_n$ must
lie between the two parallel lines (depicted in Figure \ref{fig.balltangent}) 
through $Y, Z$, while none of them will lie beneath the line through $Z$.
It follows that an edge of $\triangle A'_n B'_n C'_n$ intersects the
line $AY$ inside $\overline{B}(X,r)$.  This contradiction completes the proof.
\end{proof}

\section{The Rational Case}\label{rationalsection}

We now turn our attention to the case of rational points as in the
statement of Conjecture \ref{rationalthm}.  We note the following
simple observation.

\begin{lem}\label{6pointlemma}
Suppose that $S = \{(0,0),(0,1),(0,2),(1,0),(1,1),(1,2)\}$ or that
$S = \{(0,0),(0,1),(0,2),(1,0),(1,-1),(1,-2)\}$. Then,
$\bigcup\nolimits_{i \geq 0} {\mathcal{T}^{i} \left( S \right)} =
\mathbb Q^2.$
\end{lem}

\begin{proof}
Iteration of $\mathcal{T}$ on both sets above gives all of $\mathbb Z^2$, 
and it is easily verified that $\mathbb Z^2$ generates all of $\mathbb Q^2$.    
\end{proof}

\begin{figure}[!htbp]
\begin{center}
\includegraphics[scale=0.45]{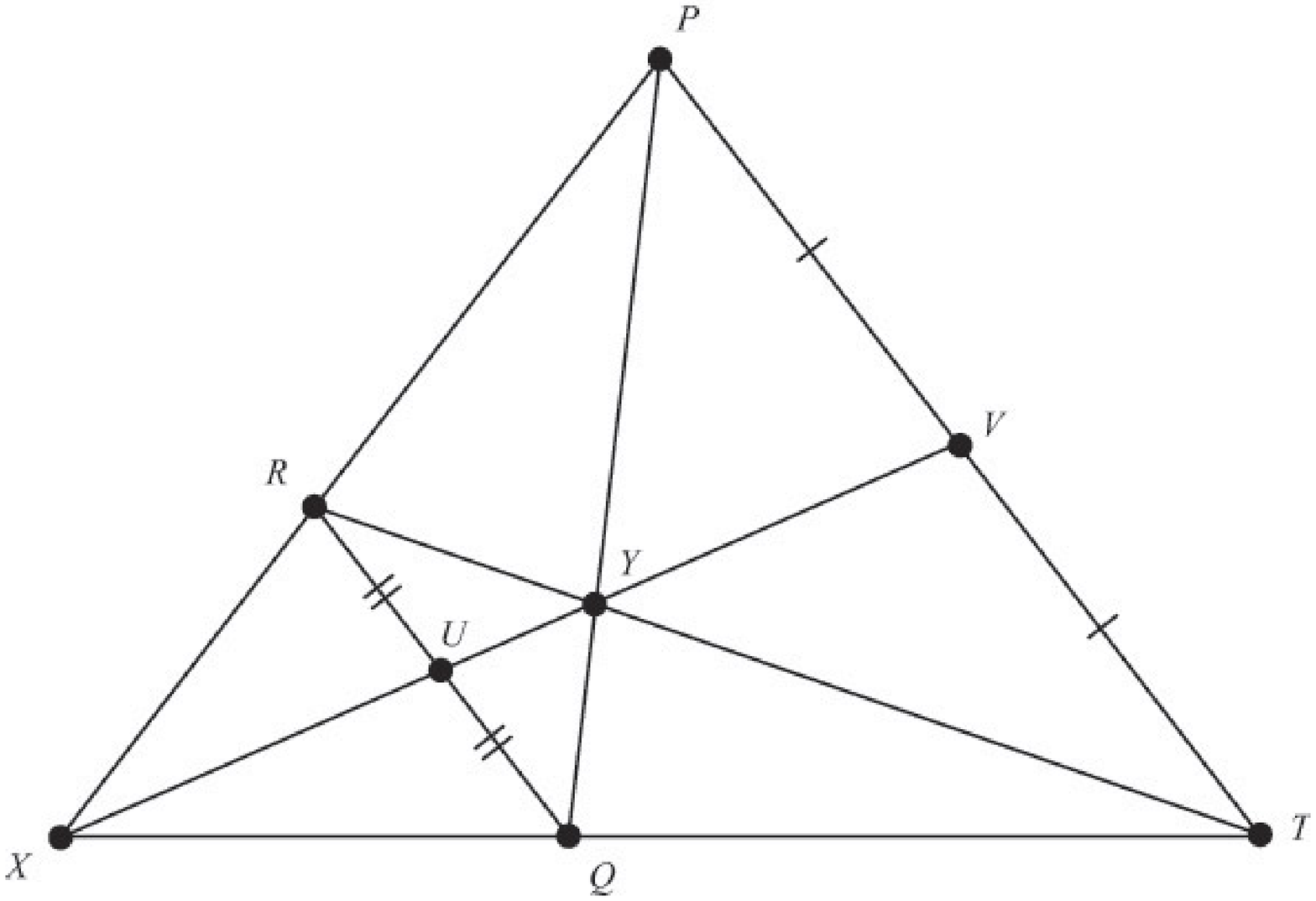}
\end{center}
\caption{Midpoint Lemma \ref{midpointlemma}} \label{fig.midpoint}
\end{figure}

We next restrict our attention to a particular case involving a
pair of parallel lines.  We need the following fact from
plane geometry.

\begin{lem}\label{midpointlemma}
Let $R, P, Q, T$ be points in the plane with $RQ$ and $PT$
parallel and suppose that $RP$ is not parallel to $QT$.  Let $Y$
be the intersection of $RT$ and $PQ$ and set $X$ to be the
intersection of $RP$ and $QT$.  Then, $XY$ intersects $RQ$ and
$PT$ in their midpoints $U$ and $V$, respectively.
\end{lem}

\begin{proof}
Since $\bigtriangleup RUY$ and $\bigtriangleup TVY$ are similar
triangles, we have $RU/TV = UY/VY$.  The same reasoning gives us
that $UY/VY = UQ/VP$.  Examining the large triangles
$\bigtriangleup XVT$ and $\bigtriangleup XPV$, it is also clear
that $UQ/TV = XU/XV = RU/VP$.  Therefore,
\begin{equation*}
UQ = TV \cdot \frac{RU}{VP} = TV^2 \cdot \frac{UQ}{VP^2},
\end{equation*}
so that $TV = VP$.  A similar computation shows that $RU = UQ$.
\end{proof}

We finally arrive at our main result in the rational case.
It will be a consequence of Lemma 
\ref{midpointlemma}, and it is the closest we come to 
proving Conjecture \ref{rationalthm}.

\begin{proof}[Proof of Theorem \ref{rationalparallelthm}]
Since a (rational) translation does not change the problem, we may
assume that $Q = (0,0)$.  Moreover, it is easy to see that if $M
\in GL_2(\mathbb Q)$, then \[M \cdot S = \left\{M \left[
{\begin{array}{*{20}c}
   a  \\
   b  \\
 \end{array} } \right] \ : \ (a,b) \in S\right\}\] gives rise to $\mathbb
Q^2$ through iteration of $\mathcal{T}$ if and only if $S$ does.
Suppose that $R = (a,b)$, $P = (c,d)$, and $T = (u,v)$ with
$a,b,c,d,u,v \in \mathbb Q$.  Since $RQ$ and $PT$ do not define
the same line, it follows that $bu-av \neq 0$.  Also, since $RQ$
and $PT$ are parallel, we have $bu-av = bc - ad$. 

Consider the following matrix:
\begin{equation*}
M = \frac{1} {{bu - av}}\left[ {\begin{array}{*{20}c}
   b & { - a}  \\
   { - v} & u  \\
 \end{array} } \right].
\end{equation*}
A straightforward computation gives $M \cdot S =
\left\{(0,0),(0,1),(1,0),\left(1,\frac{du-cv}{bu-av}\right)\right\}.$
Moreover, since $RP$ is not parallel to $QT$, it follows that
$\frac{du-cv}{bu-av} \neq 1$.  Next, set $\frac{r}{s} =
\frac{du-cv}{bu-av}$ in which $r,s \in \mathbb Z$ and gcd$(r,s) =
1$. Suppose first that $r/s > 0$.  By successively applying Lemma
\ref{midpointlemma}, iteration of $\mathcal{T}$ on $M \cdot S$ produces
the points:
\begin{equation*}
\left\{\left(0,\frac{l_1}{2^k}\right),\left(1,\frac{rl_2}{s2^k}\right)
\ : \ l_1,l_2,k \in \mathbb N; \ 0 \leq l_1,l_2 \leq 2^k \right\}.
\end{equation*}
It follows that if we choose $k$ such that $2^{k-1} \geq
\text{max}\{r,s\}$, we will have
\begin{equation*}
\left\{\left(0,0\right),\left(0,\frac{r}{2^k}\right),\left(0,\frac{2r}{2^k}\right),
\left(1,0\right),\left(1,\frac{r}{2^k}\right),\left(1,\frac{2r}{2^k}\right)\right\}
\subseteq \mathcal{T}^{k} \left( M \cdot S \right).
\end{equation*}
Therefore, letting $N = \left[ {\begin{array}{*{20}c}
   1 & 0  \\
   0 & {\frac{{2^k }}
{r}}  \\
 \end{array} } \right]$, we must have
\begin{equation*}
\left\{(0,0),(0,1),(0,2),(1,0),(1,1),(1,2)\right\} \subseteq N
\cdot \mathcal{T}^{k} \left( M \cdot S \right).
\end{equation*}
An application of Lemma \ref{6pointlemma} now concludes the proof
of this case.  

Finally, if $r/s < 0$, then the same examination as above reduces
the situation to $S = \{(0,0),(0,1),(0,2),(1,0),(1,-1),(1,-2)\}$, also covered by Lemma \ref{6pointlemma}.
\end{proof}

\section{Acknowledgments}

We would like to thank the anonymous referees for introducing us to the references 
\cite{barany, pach, otherproof, king} and for several suggestions that 
improved the exposition of this paper.  Special thanks also go to Kelli Carlson
for helping to simplify the proof of Lemma \ref{trianglelemma2} and for
giving useful comments on preliminary versions of this paper.



\begin{thebibliography}{1}

\bibitem{bookproof}
M. Aigner and G. M. Ziegler, \emph{Proofs from the book},
Springer-Verlag, New York, 1999.

\bibitem{barany}
I. Barany, P. Frankl, H. Maehara, \emph{Reflecting a triangle in the plane}, Graphs and Combinatorics,
\textbf{9} (2), 97--104, 1993.

\bibitem{pach}
K. Bezdek and J. Pach, \emph{A point set everywhere dense in the plane},
Elem. Math., \textbf{40} (4), 81-84, 1985.

\bibitem{solvedsylvester}
T. Grunwald, \emph{Solution to Problem 4065}, Amer. Math. Monthly, \textbf{51}, 169-171, 1944.

\bibitem{otherproof}
D. Ismailescu and R. Radoicic, \emph{A dense planar point set from iterated 
line intersections}, Comp. Geom. Theor. Appl., \textbf{27} (3), 257--267, 2004.

\bibitem{king}
J. King, \emph{Three problems in search of a measure}, 
Amer. Math. Monthly, \textbf{101} (7), 609--628, 1994.

\bibitem{numtheory}
I. Niven, H. Zuckerman, H. Montgomery, \emph{An introduction to
the theory of numbers}, John Wiley \& Sons, Inc., New York, 1991.

\end{thebibliography}
\end{document}